# Geometry on Optimal Problem

Denys Shcherbak

denys.shcherbak@umu.se

Natalya Pya Arnqvist

natalya.pya@umu.se

Department of Mathematics and Mathematical Statistics,

Umeå University, 901 87 Umeå, Sweden

December 4, 2023

**Abstract**

We introduce an algorithm which can be directly used to feasible and optimum search in linear programming. Starting from an initial point the algorithm iteratively moves a point in a direction to resolve the violated constraints. At the same time, it ensures that previously fulfilled constraints are not breached during this process. The method is based on geometrical properties of $n$-dimensional space and can be used on any type of linear constraints $(>, =, \geq)$, moreover it can be used when the feasible region is non-full-dimensional.

## 1 Introduction

Consider a system of linear inequalities

$$\begin{cases} a_{11}x_1 + a_{12}x_2 + \cdots + a_{1n}x_n & \geq \quad b_1 \\ a_{21}x_1 + a_{22}x_2 + \cdots + a_{2n}x_n & \geq \quad b_2 \\ \vdots & \quad \vdots \\ a_{m1}x_1 + a_{m2}x_2 + \cdots + a_{mn}x_n & \geq \quad b_m \end{cases} \qquad (1)$$

with an objective function

$$f(x) = c_{11}x_1 + c_{12}x_2 + \cdots + c_{1n}x_n \qquad (2)$$

One of the questions that might arise here, is whether any solution of the system of inequalities exists? If yes, how to minimize the value of $f(x)$ over all possible solutions. This is a classical formulation of *Linear Programming* problem (LP). The inequalities in (1) are called *Constraints* and $f(x)$ is called *Objective function*. The system (1 )can be written in a matrix form as follows:

$$Ax \geq b$$

where $a_i = (a_{i1}, a_{i2}, \ldots, a_{in})$, $c = (c_{i1}, c_{i2}, \ldots, c_{in})$ and $b = (b_1, b_2, \ldots, b_m)$.

The starting point of LP is considered to be the report [1] by young Soviet professor Leonid Kantorovich in 1939. The report was on organization and production planing at plywood trust laboratory, where Kantorovich formulated an extremum problem under a system of linear inequalities. An interesting fact is that, in 1930-40s the only way to work with linear inequalities was Fourier-Motzkin elimination method. The idea of eliminating variables from a system of linear inequalities was first introduced by Fourier in 1826 [2]. This concept was independently rediscovered by Dines in 1919 [3] and then again by Motzkin in his 1936 PhD thesis [4]. Several other researchers arrived at the same idea, for more history see [5]. While this method is intuitive, it has strong restrictions due to its effectiveness. The Achilles' heel of the Fourier-Motzkin elimination is its double exponential complexity, which makes it useless for LP with its tremendous number of constraints. For example, to solve plywood trust problem described in [1] on 8 types of peeling machines and on 5 different materials it must be solved more than billion linear inequalities, which, definitely, was impossible on that time and still is a monumental challenge nowadays. For solving this problem, Kantorovich developed a method of resolving multipliers as presented in [1]. The monograph stirred significant interest in the West scientific community and was translated to English in 1960 [6], primary by the initiative of Koopmans, who worked on the theory of transportation. In recognition of their groundbreaking contributions to the theory of optimum allocation of resources, Koopmans and Kantorovich were jointly awarded the Nobel Prize in Economic Sciences in 1975.

Perhaps, the most famous method to solve LP is a simplex method, which was discovered by Dantzig in 1947 during his work on the transportation problem. The first publication of this method was in 1951 for a general case [7] and for the transportation problem [8]. The method is based on the algebraic properties of the matrix $A$ and it is quite intuitive. However it is far from simple when considering its complexity. In 1970 Klee and Minty in [9] demonstrated an example where the simplex method require an exponential number of pivoting steps. On the other hand, rapidly developed semiconductors industry and growing usage of computing machines allowed to reach new results in many areas using the simplex method, without caring about its theoretical complexity. Moreover, in practice, the simplex method demonstrated good results, far from exponential time. Thus, it is still used in almost unchanged form.

However, the simplex method was not only the way to solve extremal problem under linear inequalities. In 1976-77 Soviet scientists Nemirovski and Shor independently devised a method based on geometrical properties of $n$-dimensional space. The method is known as ellipsoid method. Both scientists came to the idea from the distinct directions. Nemirovski [10] derived the method from a central section scheme and called it as modified centred cut. Shor, in turn, was investigating a special case of his space dilation [11]. In both cases the method did not expected to be used in LP. The situation has changed dramatically in 1979, when Khachiyan presented "A polynomial algorithm in linear programming" [12]. In the paper he adapted the ellipsoid method to answer the question if a system of linear inequalities is feasible. The main result is that the answer could be given in polynomial number of steps. This, in turn, rises a natural question, whether the decision problem on system of linear inequalities is not $NP-$complete or $P = NP$.

That was the beginning of the second life of the ellipsoid method. According to [13] survey, in November 1979, The New York Times exclaimed "Shazam" in the context of Khachiyan's results. See [14] for more about the tale and the history of convex optimization. However, the ellipsoid method has one big disadvantage in practice, despite the fact that it is polynomial, it is very slow, much slower than the simplex method. The striking fact is that the non-polynomial simplex method in practice is much faster than the polynomial ellipsoid method.

In this work we want to present a method on system of linear inequalities which is based on geometrical properties of $n$-dimensional space. In Section 2 we show simple examples how the method could work in 2 and 3-dimensional spaces, where human imagination allows us to understand the behaviour and the logic of the method. The examples are straightforward and purposed to give an intuition of how the method acts in $n$-dimensional space. The main algorithm is described in Section 3. Furthermore, in this section we show some narrow places of the algorithm and how to resolve the corresponding issues in theoretical and practical sides. The main outline of the of Sections 2 and 3 is to resolve violated constraints, which is usually called as *Phase I* in an optimization problem. In Section 4 we show how to adapt the algorithm for optimum search and consider some practical issues of the algorithm. Finally, the concluding remarks are given in Section 5.

## 2   Preliminaries and Notations

We consider the system (1) from a geometrical prospective. In other words, we look at $Ax \geq b$ as a body in $n$-dimensional space. Its convex hull is defined by hyperplanes $a_{i1}x_1 + a_{i2}x_2 + \cdots + a_{in}x_n = b_i$ for $i = 1, 2, \ldots m$. We denote by $l_i$ the corresponding hyperplane, for simplicity we will refer to it as a plane. It is easy to see, that the $i$-th constraint defines a half-space. The normal vector to $i$-th plane we denote by $\vec{a_i} = (a_{i1}, a_{i2}, \ldots, a_{in})$. Without loss of generalisation we may assume that all vectors $\vec{a_i}$ are normalized, if not, we can update the system of inequalities dividing each constraint by corresponding length of vector $||\vec{a_i}||$.

Historically, the word **vector** has several contexts. From an algebraic point of view vectors $a_i$ and $b$ are the same objects, namely all operations on them are defined identically, the difference is only the dimensions. However, vector $a_i$ has a geometric sense: it indicates the direction. But $b$ has no such property. To distinguish these **algebraic** and **geometric** vectors we use arrow $\vec{}$ to emphasize that vector describes a certain direction and we use its geometric nature.

To answer the question about consistency of (1), one could slightly modify the original problem and use the simplex method or use Khachiyan's variant of the ellipsoid method. In both cases, an algorithm will construct a sequence of points which will lead to a feasible region.

In the following sections we describe the algorithm, which aims to find a feasible solution. Starting with an initial point, algorithm constructs a sequence of points, which safely reduce the number of violated constraints at each iteration, and reach feasible region. We start with the simplest case and eventually increase the complexity and finally, introduce the general description of the method in Section 3.

## 2.1 The Simplest Case: 1-constraint

Assume $m = 1$, meaning that there is a single constraint, and we need to find a feasible solution. Let $P_0$ be a point. If it satisfies the constraint, we are done – we found the point in the feasible region.

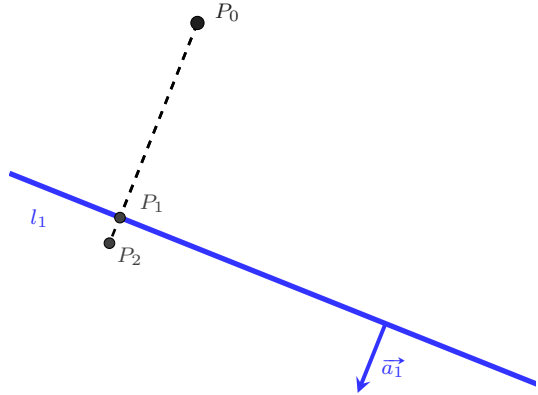

Figure 1: Single constraint.

Let $P_0$ violates the constraint. Fig. 1 schematically shows the corresponding plane $l_1$. The normal vector $\vec{a}1$ shows the direction where the constraint is fulfilled. In this simple case we can easily determine a point which satisfies the constraint. Obviously, we can find a projection of $P_0$ on the plane $l_1$, namely point $P_1$ in Fig 1. If the inequality is not strict ($\geq$), $P_1$ satisfies the constraint. If not (namely $>$), we can move $P_1$ a bit further along the vector $\vec{a_1}$ on distance $\epsilon$ and find point $P_2 = P_1 + \epsilon \cdot \vec{a_1}$. This process has a trace $P_0 \to P_1 \to P_2$, that we will naturally refer to as a movement: "move point $P_0$ to $P_2$, which satisfies the constraint".

## 2.2 A Simple Case: 2-constraints

Consider a system of two inequalities, is it consistent or not? In other words, we need to find a point (if any) which satisfies both inequalities. Fig. 2 schematically shows two inequalities: theirs planes $l_1, l_2$ and normal vectors $\vec{a_1}, \vec{a_2}$, which show the direction where the corresponding inequalities are fulfilled. It is easy to see, that this system is consistent, the feasible region is coloured in green.

If the number of variables were $n = 2$, there is no need to worry. In such a case we could just find an intersection of $l_1$ and $l_2$, and then move that point along the vector $\vec{w} = \vec{a_1} + \vec{a_2}$. This will directly lead the point being in the feasible region. However, the problem is more difficult in general case, as the intersection of two planes in $n$-dimensional space is not a single point.

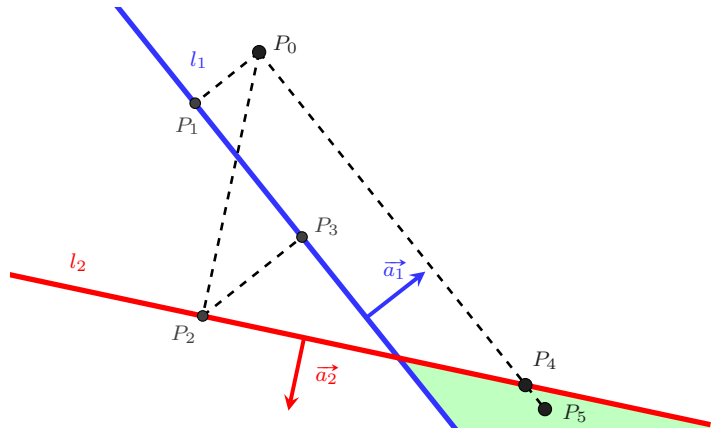

Figure 2: Couple of constraints.

Let $P_0$ be an initial point. If it satisfies both constraints, the problem is solved as $P_0$ is in feasible region. If $P_0$ violates both constraints, as it has been shown before we can easily move the point such that at least one of the constraints is fulfilled. Thus, let us assume that $P_0$ satisfies the first constraint and violates the second one as in Fig 2. To resolve the second constraint one may move $P_0$ to point $P_2$. In such a case the movement will cross plane $l_1$, and consequently, the first constraint will be violated. The aim of our method is to resolve the second constraint while keeping the first constraint fulfilled. In other words, we would move the point such that it will not cross the plane $l_1$. Obviously, to omit the risk of crossing, the best direction of the movement is parallel to $l_1$. Perhaps, the easiest way to construct the vector of such direction is to use normals $\vec{a_1}$, $\vec{a_2}$ and $\cos(\widehat{\vec{a_1}, \vec{a_2}})$. However, we present another way to construct it. It is not very natural for two constraints, but we use similar idea in the following sections. To construct the vector of such movement, we need to compute the following:

- find $P_1$ as a projection of $P_0$ on $l_1$,
- find $P_2$ as a projection of $P_0$ on $l_2$,

- find $P_3$ as a projection of $P_2$ on $l_1$.

It is easy to see that the vector $\overrightarrow{P_1P_3}$ is indeed the desired direction. Having point $P_0$ and the direction vector $\overrightarrow{P_1P_3}$ we can define a line. It is not a problem to find the intersection point $P_4$ of the defined line and the plane $l_2$. Similarly, as have we did earlier, if constraint 2 is not strict, then $P_4$ satisfies both constraints. If not, we can move $P_4$ a bit further by $\epsilon$ and get $P_5 = P_4 + \epsilon \frac{\overrightarrow{P_1P_3}}{||\overrightarrow{P_1P_3}||}$, where $\epsilon$ is a small value.

## 2.3 A Hard Case: 3-constraints

Consider now the system with three constraints. Is the system consistent? Let $P_0$ be an initial point. By repeating the similar arguments as in the previous sections, we can reach the situation when $P_0$ satisfies at least two of the constraints. Let's say, $P_0$ satisfies constraints 1 and 2, but violates the constraint 3. At Fig. 3 schematically shown the corresponding planes $l_1, l_2, l_3$, their normals $\vec{a_1}, \vec{a_2}, \vec{a_3}$ and the initial point $P_0$.

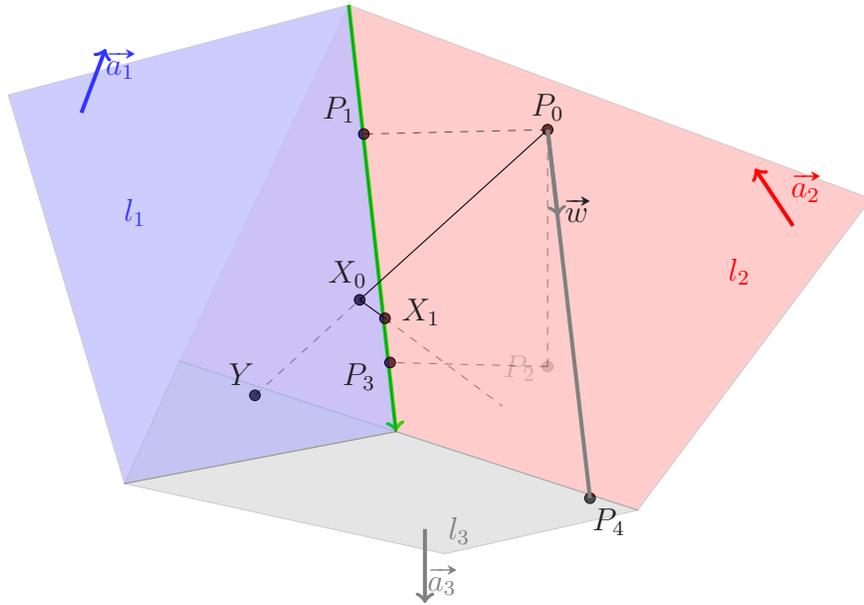

Figure 3: Schematic of 3 constraints.

It would be nice if we could repeat the procedure for two constraints as in the previous section to resolve the third constraint. However, it may happen that by moving along one of the planes we may cross another one. For example, if we move $P_0$ parallel to plane $l_2$ with an aim to cross plane $l_3$ (line $P_0Y$, $Y \in l_3$ Fig. 3) we will cross plane $l_1$ at point $X_0$. Obviously, $Y$ is not a solution for the system. On the other hand, we could stop the movement at point $X_0$ on plane $l_1$. In such a case, $X_0$ still satisfies constraints 1 and 2, but violates constraint 3. If we would repeat this procedure with $X_0$, we will cross plane $l_2$ at point $X_1$, and this is a stack.

As previously, we would move $P_0$ such that it will not crosses the planes, namely parallel to planes $l_1, l_2$. The only question, is how to do it simultaneously to both planes. In 3-dimensional space, it could be performed effortlessly by finding a vector product of the normal vectors $\vec{w} = \vec{a_1} \times \vec{a_2}$. The vector $\vec{w}$ guarantees that such a movement will not cross either $l_1$ or $l_2$.

However, there are two issues here. Firstly, strictly speaking, the intersection of two planes in $n$-dimensional space is not a line. In addition, the vector product is not defined in $n$-dimensional space. So, how to find the direction of $P_0$ which does not cross $l_1, l_2$?

For this purpose we need to find 2 points on the intersection of planes $l_1$ and $l_2$. This could be done very similar to 2-constraint case considered above. For that we need some extra computations:

- find a point $P_1$, which is the closest to $P_0$ on the intersection of $l_1$ and $l_2$,

- find a point $P_2$ as a projection of $P_0$ on $l_3$,

- find a point $P_3$, which is the closest to $P_2$ on the intersection of $l_1$ and $l_2$.

The vector $\overrightarrow{P_1P_3}$ is indeed the direction. By the construction, $\overrightarrow{P_1P_3}$ is parallel to $l_1$ and $l_2$, thus it does not cross either of them. Moreover, the point $P_3$ is closer to $l_3$ than $P_1$, thus $\overrightarrow{P_1P_3}$ crosses $l_3$. Thus, we can move $P_0$ along the vector $\overrightarrow{P_1P_3}$ and find a point $P_4$ on the plane $l_3$. It is clear that $P_4$ fulfils all three constraints. Note, it might happen, that $P_3$ not closer to $l_3$, namely $P_3 = P_1$, which means that the system is not feasible.

**Observation 2.1.** *Actually we can reach the feasible region without simultaneous parallel moving. As we show, first we can move along $l_1$. Once we intersect $l_2$ we can move parallel to it. By continuing proposed "zig-zag" procedure we will reach plane $l_3$ and resolve the constraint. However, there is an* *important restriction: the "zig-zag" procedure will work if the points of the sequence $X_0, X_1 \ldots$ are inner points, for both constraints 1 and 2. In other words, the points $X_0, X_1 \ldots$ should not be on the planes, but instead above them. To achieve that, we can make a small $\epsilon$-shift along the corresponding normal vector. It is clear, that the procedure is require a lot of such "zig-zag" steps. Taking in account that all calculations are done in n-dimensional space, this method is not very attractive due to its computational load.*

## 2.4 The Closest Point

In Section 2.3 we defined points $P_1, P_3$ as the closest points to $P_0, P_2$ respectively, which lie on the intersection of planes $l_1, l_2$. Note, that in this case $P_1, P_3$ are projections of $P_0, P_2$ on the intersection. To determine the coordinates of $P_1, P_3$ we will use a pure algebraic technique: a least-square solution method.

Let $l_1, l_2, \ldots, l_t$ be planes with non empty intersection. Furthermore, let $P = (p_1, p_2, \ldots, p_n)$ be a point. Let us find a point $T$ on the intersection of the planes and which is the closest possible to $P$. The plane $l_i$ is defined by equation $a_{i1}x_1 + a_{i2}x_2 + \cdots + a_{in}x_n - b_i = 0$. Denote by $G$ a matrix of normal vectors of the planes

$$G := \begin{bmatrix} a_1 \\ a_2 \\ \vdots \\ a_t \end{bmatrix}.$$

Denote as $d_i(P) := a_{i1}p_1 + a_{i2}p_2 + \cdots + a_{in}p_n - b_i$. The coordinates of the point $T$ can be found as a least-square solution:

$$T = G^T(GG^T)^{-1}d =$$

$$= \begin{bmatrix} | & | & \ldots & | \\ a_1 & a_2 & \ldots & a_t \\ | & | & \ldots & | \end{bmatrix} \cdot \begin{bmatrix} \vec{a_1} \cdot \vec{a_1} & \vec{a_1} \cdot \vec{a_2} & \ldots & \vec{a_1} \cdot \vec{a_t} \\ \vec{a_2} \cdot \vec{a_1} & \vec{a_2} \cdot \vec{a_2} & \ldots & \vec{a_2} \cdot \vec{a_t} \\ \vdots & \vdots & \ldots & \vdots \\ \vec{a_t} \cdot \vec{a_1} & \vec{a_t} \cdot \vec{a_2} & \ldots & \vec{a_t} \cdot \vec{a_t} \end{bmatrix}^{-1} \cdot \begin{bmatrix} d_1(P) \\ d_2(P) \\ \vdots \\ d_t(P) \end{bmatrix}$$

where $\vec{a_i} \cdot \vec{a_j} = a_{i1}a_{j1} + a_{i2}a_{j2} + \cdots + a_{in}a_{jn}$ is a scalar product of two vectors.

# 3 Direction Search Algorithm

Now we can generalize the ideas outlined in the previous sections. Let $Ax \geq b$ be a system of $m$ linear inequalities, and we need to find a point $P \in \mathbb{R}^n$, which satisfies the system. Similarly to the previous examples, we will move the point to resolve the violated constraints. However, unlike the sections 2.1 - 2.3, we consider the point $P$ with its $\epsilon$-neighbourhood, namely $\epsilon$-ball with a centre at $P$. The reason of such a restriction and the ways how to omit it we will discuss in later sections.

Denote as $d_i(P) := a_{i1}p_1 + a_{i2}p_2 + \cdots + a_{in}p_n - b_i$. Since all $\vec{a_i}$ are normalized, $d_i(P)$ is almost classical (metric) distance between a point $P$ and a plane $l_i$. The only difference is that the classical distance is always non-negative. However, sometimes it is useful to know the sign of the value. For example, if the distance is negative $d_i(P) < 0$ this means that $P$ violates $i$-th constraint. For simplicity we will call it also a distance. We denote by $r(P, \overrightarrow{dir})$ a ray with a starting point $P$ which is parallel to a vector $\overrightarrow{dir}$.

The main idea of the algorithm is (gradually) to resolve the violated constraints and not braking the fulfilled constraints. In other words, during the movement of point $P$ we allow to intersect only the planes with (strictly) negative distance to $P$, $d_i(P) < 0$.

Similar to rainwater flows down a gutter, moving along its walls, our $\epsilon$-ball will move along to planes in the $n$-dimensional space. The only difference is that the number of planes which affect the $\epsilon$-ball might be more than three. Naturally, the set of such planes we will call as a ***gutter***.

Let $P_0$ violate the $i$-th constraint, which means $d_i(P_0) < 0$. To describe the proposed algorithm we will use the following notations:

$P_0$ is the centre of $\epsilon$-ball;

$P_2$ is the projection of $P_0$ on $l_i$ (violated constraint);

$\overrightarrow{dir}$ is a vector of movement;

$r(P_0, \overrightarrow{dir})$ is the ray starting at $P_0$ in the direction $\overrightarrow{dir}$;

$G$ is a matrix of normal vectors of planes which form a gutter;

$u$ is a vector of distances between $P_0$ and planes which form a gutter (all $\epsilon$'s);

$v$ is a vector of distances between $P_2$ and planes which form a gutter.

The dimensions of $G, u$ and $v$ will depend on how many planes affect the movement of the $\epsilon$-ball at current iteration, and could be different at each iteration step. How to define the elements of the matrix and vectors will be described in the algorithm below. Note, $u$ and $v$ are algebraic vectors and their dimensions are not fixed, however $\overrightarrow{dir}$ is a geometric vector which determines a direction of a movement.

We will now present the sketch of $i$-th constraint resolving process, which generalizes the ideas from the previous sections. The formal description of the direction search algorithm is shown as Algorithm 1 below.

*Initial step:* There is no gutter, no plane which affects the direction of $\epsilon$-ball. Thus, set $\overrightarrow{dir} = \vec{a_i}$ and set $G, u, v$ to be empty.

*First step:* Move $\epsilon$-ball in direction $\overrightarrow{dir}$ as close as possible to $l_i$. By "as close as possible" we mean not braking the constraints which $P_0$ already fulfils. After this step, either we resolve $i$-th constraint, or $\epsilon$-ball touches a plane from a gutter.

*Gutter step:* When $\epsilon$-ball ends up at one of the planes (which forms a gutter) we do updates of element of $G$, $u$ and $v$ as follows:

- append vector $a_j$ to $G$;
- append $d_j(P_0) = \epsilon$ value to $u$;
- append the $d_j(P_2)$ to $v$.

Then, we update the coordinates of $P_1, P_2, P_3$, where $P_1$ is a projection of $P_0$ on the intersection of planes which form the gutter
$$P_1 := G^T(GG^T)^{-1}u,$$

$P_3$ is a projection of $P_2$ on the intersection of planes which form the gutter
$$P_3 := G^T(GG^T)^{-1}v.$$

Note, that the distance between $P_2$ and $l_j$ is always non-positive.

By repeating this process we will either resolve all violated constraints and find a feasible point or find an unresolvable constraint.

Since both $P_1$ and $P_3$ are on the intersection of the planes, the moving of the $\epsilon$-ball along the vector $\overrightarrow{P_1P_3}$ will not cross any of the planes which are already included to the gutter $G$. It is easy to see, that if there is a slope on the gutter, then $P_3$ will be closer to the violated plane $l_i$ than $P_1$. Thus, this movement will make $P_0$ closer to the plane $l_i$. Obviously, if there is no slope on the gutter, then $P_3 = P_1$. In such a case we would conclude that the $i$-th constraint is not resolvable and the system is infeasible.

Basing on the algebraic properties of the product $(GG^T)$ one might guess how many planes may form a gutter. Below, we provide a geometric confirmation of this fact using properties of $n$-dimensional space.

**Observation 3.1** (An upper bound of number of planes forming a gutter). *Since $\epsilon$-ball moves in a direction parallel to **one of the lines** in the intersection of the planes which form a gutter, to have at least one line in this intersection, the number of planes cannot exceed $n - 1$.*

Below we present the formal resolving algorithm.

**Algorithm 1:** Resolving violated constraints

**Data:**
Points $P_0, P_1, P_2, P_3$,
$G$ matrix on normal vectors,
$u, v$ vectors of distances,
$S, N$ sets of indices.

1 **do**
2     $N = \{i : d_i(P_0) < 0\}$     // set of violated constraints
3     $S = \{j : d_j(P_0) \geq 0\}$     // set of fulfilled constraints
4     chose $i \in N$ and put it to $S$
5     $\overrightarrow{dir} = \vec{a_i}$     // set initial direction
6     empty $G, u, v$     // resize the matrix and the vectors to 0
7     **while** *($d_i(P_0) < 0$ & number of rows in $G < n$)* **do**
8        determine the closest (to $P_0$) intersection point $X_j$ of ray $r(P_0, \overrightarrow{dir})$ and plane $l_j$, for $j \in S$
9        **if** *($j = i$)*     // the $i$-th constraint is resolved
10        **then**
11           update $P_0 = X_j + \frac{\epsilon}{\cos(\widehat{\vec{a_i}, \overrightarrow{dir}})} \cdot \vec{a_i}$     // $d_i(P_0) = \epsilon$
12           **go to** *line 1*
13        **if** *($j \neq i$)*     // obstacle on the way
14        **then**
15           **if** *($d_j(P_0) > \epsilon$)* **then**
16              update $P_0 = X_j - \frac{\epsilon}{\cos(\widehat{\vec{a_j}, \overrightarrow{dir}})} \cdot \vec{a_j}$     // $d_j(P_0) = \epsilon$
17              **go to** *line 5*
18           **else if** *($d_j(P_0) \leq \epsilon$)*     // plane of a gutter
19           **then**
20              if $d_j(P_0) < \epsilon$ update $P_0 = X_j - \frac{\epsilon}{\cos(\widehat{\vec{a_j}, \overrightarrow{dir}})} \cdot \vec{a_j}$ to get $d_j(P_0) = \epsilon$
21              detemine the coordinates of $P_2$, projection of $P_0$ on $l_i$;
22              append $d_j(P_0) = \epsilon$ to $u$;
23              append $d_i(P_2)$ to $v$;
24              append $a_i$ to $G$;
25              update the coordinates of $P_1 = G^T(GG^T)^{-1}u$;
26              update the coordinates of $P_3 = G^T(GG^T)^{-1}v$;
27              update the direction $\overrightarrow{dir} = \overrightarrow{P_1P_3}$.
28     **if** *(number of rows in $G = n$)* **then**
29        **return** *INFEASIBLE*
30 **while** *($P_0$ violates a constraint)*

## 3.1 A Non-Full-Dimensional Case

Obviously, the size of the $\epsilon$-ball should be small enough to fit the feasible region and to omit the situation illustrated in Fig. 4. Note, that for simplicity, the gutter is depicted by only two planes $l_1, l_2$. The violated constraint is presented as plane $l_3$.

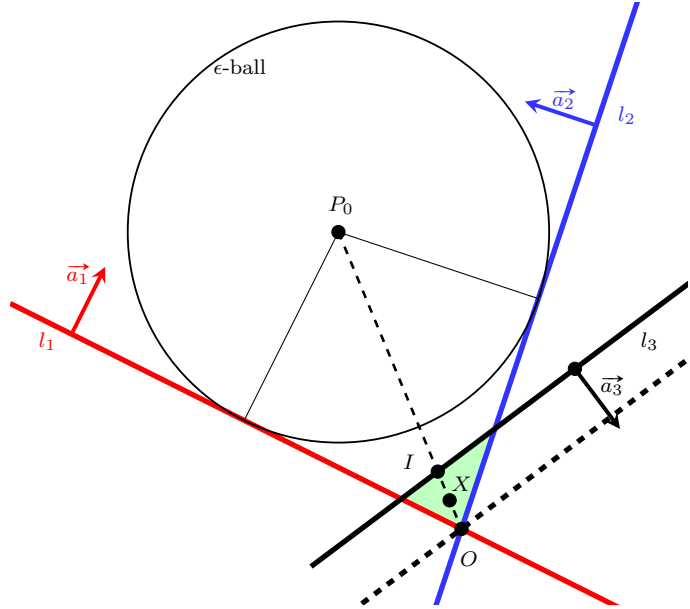

Figure 4: Large $\epsilon$-ball.

Sometimes, it is hard to estimate how small the ball should be. However, the Algorithm 1 might be adapted to resolve this issue. For this, we need to add an extra straightforward step that checks the size of the $\epsilon$-ball. If point $P_0$ satisfies the constraints of a gutter $(l_1, l_2)$ but has negative distance to $l_3$ (meaning violates $l_3$), we check if $O$ fulfils all of the constraints (see Fig. 4), where $O$ is a projection of $P_0$ on the intersection of the planes forming the gutter. The coordinates of $O$ can be found as $O = G^T(GG^T)^{-1}u$, where $u$ is a vector of all $\epsilon$'s. If $O$ does satisfy all the constraints, we can put the centre of the new $\epsilon$-ball at the point $X$, which is the middle of the segment $OI$, where $I$ is the intersection of segment $OP_0$ and plane $l_3$. The corresponding update of the radius is $\epsilon = d_1(X) = d_2(X)$, guarantees that new $\epsilon$-ball will be inside the "green zone".

However, there might be an extreme case when $l_3$ passes the intersection point $O$, see the dashed version of $l_3$ in Fig. 4, which means that the feasible region is non-full-dimensional. We can use the fact that the feasible solution is on the intersection of the gutter and the plane $l_3$, to continue the algorithm. Namely, instead of inequality ($\geq$) in the corresponding constraints during the Algorithm 1 we use equality (=). In other words, we know that the corresponding distances are equal to 0: $d_1(P_0) = 0, d_2(P_0) = 0, d_3(P_0) = 0$.

## 3.2 Why Do We Need $\epsilon$-ball?

Using $\epsilon$-ball imposes a strong restriction that the convex body $Ax \geq b$ must be full-dimensional. There is a temptation of not using the $\epsilon$-ball in the process. In other words, instead of moving the point $P_0$ close enough to the planes, move it directly on the planes. Since we do not cross the planes, Why not? Below we show a counter-example when the Algorithm 1 will fail without $\epsilon$-shift away from a plane.

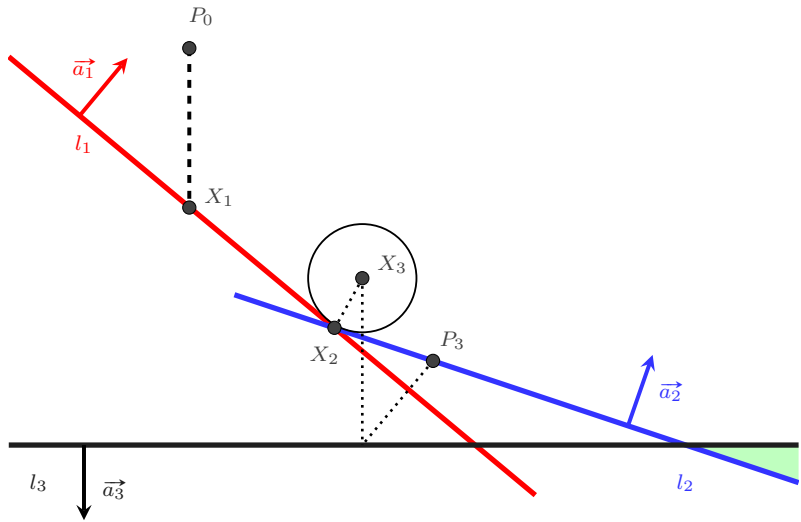

Figure 5: A counter-example on three constraints.

Let $P_0$ be an initial point which satisfies constraints 1 and 2 but violates constraint 3 (see Fig. 5). Let us apply the Algorithm 1 without $\epsilon$-shift.

Following the Algorithm 1, we move the initial point $P_0$ to $X_1$ on the plane $l_1$. The next iteration will move $X_1$ along $l_1$ until it crosses plane $l_2$ at the point $X_2$. It is easy to see that to reach the feasible region (depicted in green Fig. 5) the point should move along the plane $l_2$. However, for Algorithm 1 this is a stack, it cannot chose the correct direction. According to Algorithm 1 we need to detect which plane is the first obstacle on the direction $\vec{dir} = \vec{a_3}$. Since $X_2$ belongs to both $l_1$ and $l_2$ planes the algorithm cannot determine which plane is closer to $X_2$. As consequence, algorithm decides that $l_1$ and $l_2$ form a gutter. Obviously, if the intersection of planes $l_1$ and $l_2$ is parallel to plane $l_3$, we cannot get closer to $l_3$. Thus, the algorithm will return that the system is infeasible, despite the fact that there is a full-dimensional feasible region.

That is why we use a small shift away from the planes, which resolves that issue. It is easy to see, that the point $X_3$ does not have such an issue (see Fig. 5), we can easily determine that the ray $r(X_3, \vec{dir})$ intersects plane $l_2$ the first, thus the ball will move parallel the plane $l_2$ and reach the feasible region. The coordinates of $X_3$ can be determined as $X_3 = X_2 + \epsilon \cdot \vec{a_1} + \epsilon \cdot \vec{a_2}$.

## 4 Application

### 4.1 Optimum Search

Applying the Algorithm 1 on a system of linear inequalities (1) will search a feasible solution, which is called as *Phase I* in an optimization problem. It is easy to see, that we can directly apply the Algorithm 1 for searching an optimal value of a function $f(x)$, so called *Phase II* in an optimization problem. For this, we update the system (1) by adding extra constraint $f(x) \geq M$, where $M$ is a "large number". To be consistent in the geometrical sense we might choose the value $M$ such that the plane $f(x) = M$ be "outside" the body $Ax \geq b$ (see Fig. 6). However, strictly speaking, the algorithm only requires the normal vector of the plane $f(x) = M$. In other words, for applying Algorithm 1 $M$ can be any.

The algorithm will try to resolve the violated constraint $f(x) \geq M$. Due to the fact that the Algorithm 1 uses $\epsilon$-ball, the movement process will stop at point $P$, which is not indeed extremum. The coordinates of optimum point computed as $O = G^T(GG^T)^{-1}u$, where $u$ is a vector of all $\epsilon$'s.

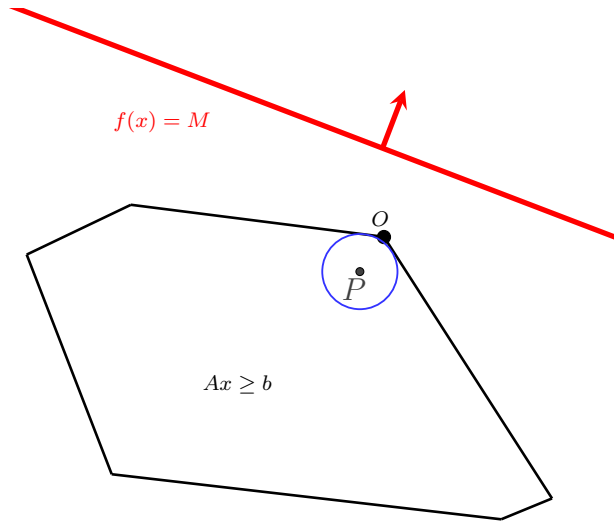

Figure 6: Optimum search.

## 4.2 Matrix Inverse

From the practical prospective the hardest stage of the algorithm is a computing of $G^T(GG^T)^{-1}$. For anyone, who are familiar with large matrices, the expression $(GG^T)^{-1}$ might trigger a panic attack. But let's look a bit details and slow the heart rate down.

First of all, $(GG^T)$ is a symmetric matrix. The second, and most important, the size of $G$ rises eventually row by row with the size of a gutter, thus so does the matrix $(GG^T)$. Use these facts to facilitate the task.

To get an inverse we can use Gaussian elimination.

$$(GG^T) = \begin{bmatrix} \vec{a_1} \cdot \vec{a_1} & \vec{a_1} \cdot \vec{a_2} & \dots & \vec{a_1} \cdot \vec{a_t} \\ \vec{a_2} \cdot \vec{a_1} & \vec{a_2} \cdot \vec{a_2} & \dots & \vec{a_2} \cdot \vec{a_t} \\ \vdots & \vdots & \dots & \vdots \\ \vec{a_t} \cdot \vec{a_1} & \vec{a_t} \cdot \vec{a_2} & \dots & \vec{a_t} \cdot \vec{a_t} \end{bmatrix}$$

$$\left[ \begin{array}{cccc|cccc} \vec{a_1} \cdot \vec{a_1} & \vec{a_1} \cdot \vec{a_2} & \dots & \vec{a_1} \cdot \vec{a_t} & 1 & 0 & \dots & 0 \\ \vec{a_2} \cdot \vec{a_1} & \vec{a_2} \cdot \vec{a_2} & \dots & \vec{a_2} \cdot \vec{a_t} & 0 & 1 & \dots & 0 \\ \vdots & \vdots & \dots & \vdots & & & \ddots & \\ \vec{a_t} \cdot \vec{a_1} & \vec{a_t} \cdot \vec{a_2} & \dots & \vec{a_t} \cdot \vec{a_t} & 0 & 0 & \dots & 1 \end{array} \right]$$

Gently applying the elimination on the first $t-1$ rows and columns we get the following:

$$\begin{bmatrix} 1 & 0 & \ldots & 0 & r_1 & & & & 0 \\ 0 & 1 & \ldots & 0 & r_2 & & & & 0 \\ \vdots & \vdots & \ddots & \vdots & \vdots & & Q & & \vdots \\ 0 & 0 & \ldots & 1 & r_{t-1} & & & & 0 \\ \vec{a_t} \cdot \vec{a_1} & \vec{a_t} \cdot \vec{a_2} & \ldots & & \vec{a_t} \cdot \vec{a_t} & 0 & 0 & \ldots & 1 \end{bmatrix} \quad (3)$$

where values $r_i$ are receiving during the computation. The matrix $Q$ is the inverse for corresponding gutter on $t-1$ planes. This means that we can use the data from the previous stages to simplify computations. In other words, there is no need compute inverse from scratch at each iteration. Once $G$ is updated by one extra row we need to fill $t$-th row by getting the values $\vec{a_t} \cdot \vec{a_i}$ for $i = 1, 2 \ldots t$, and fill $t$-th column by getting the values

$$r_i = Q(\vec{a_t} \cdot \vec{a_1}, \vec{a_t} \cdot \vec{a_2}, \ldots, \vec{a_t} \cdot \vec{a_{t-1}})^T$$

for $i = 1, 2 \ldots t - 1$. To get the inverse we need to eliminate only the $t$-th row and the $t$-th column in 3, but not the whole $t \times t$ matrix.

## 5 Conclusion

In this work we presented the algorithm, which can be directly applied to both phases of optimisation problem. It is clear, that the authors attempting to generalize a ball movement in $n$-dimensional space under certain restrictions. Due to the luck of knowledge in the field of theoretical mechanics, and particularly its extension to $n$-dimensional space, we have not taken in account any transitional process, such as inertia. In this regard, we have a couple of open questions, which might improve the implementation of the desired algorithm. At the Fig. 5 depicted a ball of radius $\epsilon$ with the centre in $X_3$. For simplicity, we may assume that $\vec{a_3}$ is a gravitation and planes $l_1, l_2$ are surfaces which are orthogonal to sheet, where $\vec{a_1}, \vec{a_2}$ represent surfaces reactions. In general vectors $\vec{a_1}, \vec{a_2}, \vec{a_3}$ could be any types of force, which affect the ball.

**Question 5.1.** *The ball at the initial moment touches both red and blue surface. Obviously, the ball will roll down the blue surface. The question is: How does (if does) the red surface affect the initial movement of the ball?*

The second question concerns an extremal case.

**Question 5.2.** *What happens if radius of the ball is zero? In other words, how will move the point $X_2$ under effect of forces $\vec{a_1}, \vec{a_2}, \vec{a_3}$?*

Our approach is based on the geometrical properties of $n$-dimensional space and draws inspiration from the ellipsoid method, and especially from Khachiyan's genius adaptation for feasible solution search. The core feature of ellipsoid method on optimisation problems is to control, in certain way, the volumes of nested ellipsoids which contain a feasible region. The scientific importance of this method is that the sequence of such ellipsoids could be constructed in polynomial time. However, despite the fact that it is theoretically fast, the squeezing of the ellipsoids occurs too slowly, and this, in turn, greatly restricts the method's applicability. Khachiyan's method mainly focuses on computation of ellipsoids in $n$-dimensional space (which take the major of computing time), but not on feasible solution or optimum value. In contrast to the ellipsoid method, our algorithm truly concentrate on feasible solution or optimum value, depending on phase of an optimisation problem.

Another state-of-art algorithm on optimisation problem is Dantzig's simplex method. The number of applications for simplex method cannot be quantified. The root of the method is algebraic properties of matrix determinant. Fast transforming matrix rows and columns is carried out with the goal of updating the objective function's value and ultimately reaching an optimal solution. As history shows the word "fast" is the key feature for computation. However, it's worth noting that the simplex method also has some restrictions. The main of them is the simplex method theoretically non-polynomial. The pivoting process is non-trivial, and it typically involves applying a greedy approach during the column pivoting step. Practical experience shows that using a non-greedy approach may impact the number of iterations, and this impact not necessary leads to an increase number of steps. In other words, the simplex method faces difficulty in determining the optimal direction, and this, perhaps, is the main contributing factor to its exponential behaviour. In a contrast to simplex method, our method does not have such a problem, it choose the unique direction at each iteration, moreover the direction is the shortest and fastest path toward achieving an optimal value. Another feature of the simplex method is variables, depends on the problem it might require extra variables. Beyond the fact that this myriad of variable types can be challenging for young researchers in learning this method, it significantly increases the size of operational matrix, which is sensitive for computation resources, especially for memory. Our algorithm does not require any extra variables it works only with the variables from origin problem. Thus, the memory storage is exactly the same as problem dimension, namely it uses $m \times n$ cells for matrix $A$, $m$ cells for vector $b$ and $n$ cells for objective function's coefficients. Moreover, our method works directly with any type of inequalities $(>, =, \geq)$ without any modifications of constraints enhancing the usability. Finally, there is no restriction on the boundedness of convex hull, this is not the main feature for optimisation problem, however it might find an application in closely related problems.

Overall, our algorithm is straightforward and intuitively clear, at the same time it allows to omit several important obstacles known in other classical methods.